\documentclass[a4paper,11pt]{article}
\usepackage{amsmath,amsfonts,amssymb,times,mathptmx,graphicx,amsthm,lineno}
\usepackage{placeins}
\usepackage[T1]{fontenc}
\usepackage[latin1]{inputenc}
\title{\bf 
 Proofs of non-optimality of the standard least-squares method for track reconstructions }
\author{Gregorio Landi$^a$\thanks{Corresponding
author. Gregorio.Landi@fi.infn.it}~,   Giovanni E. Landi$^b$\\
\\
\llap{$^a$} Dipartimento di Fisica e Astronomia,
Universita' di Firenze and INFN\\
Largo E. Fermi 2 (Arcetri) 50125, Firenze, Italy\\
\\
\llap{$^b$} ArchonVR S.a.g.l.,\\
Via Cisieri 3,
6900 Lugano, Switzerland.\\ \\
}
\date{ March 22, 2020 }
\begin{document}
\maketitle 
\begin{abstract}
It is a standard criterium in statistics to define an optimal
estimator the one with the minimum variance. Thus, the optimality
is proved with inequality among variances of competing estimators.
The inequalities, demonstrated here, disfavor the standard least squares
estimators.
Inequalities among estimators are connected to names of
Cramer, Rao and Frechet. The standard demonstrations of these inequalities
require very special analytical properties of the probability functions, globally
indicated as regular models. These limiting conditions are too restrictive
to handle realistic problems in track fitting. A previous extension to
heteroscedastic models of the Cramer-Rao-Frechet inequalities was
performed with Gaussian distributions. These demonstrations proved
beyond any possible doubts the superiority of the heteroscedastic models
compared to the standard least squares method. However, the
Gaussian distributions are typical members of the required regular
models. Instead, the realistic probability distributions, encountered
in tracker detectors, are very different from Gaussian
distributions. Therefore, to have well grounded set of inequalities,
the limitations to regular models must be overtaken. The aim of
this paper is to demonstrate the inequalities for
least squares estimators for irregular models of probabilities,
explicitly excluded by the Cramer-Rao-Frechet demonstrations.
Estimators for straight and parabolic tracks will be considered.
The final part deals with the form of the distributions of simplified
heteroscedastic track models reconstructed with optimal estimators
and the standard (non-optimal) estimators. A comparison among
the distributions of these different estimators shows the large loss in
resolution of the standard least squares estimators.

\end{abstract}



\newpage
\tableofcontents

\pagenumbering{arabic} \oddsidemargin 0cm  \evensidemargin 0cm


\section{Introduction}\indent

It is evident the crucial importance to study  the inequalities
among estimators of different fitting algorithms.
In fact, a standard definition of optimality is connected to the
estimator with the minimum variance among competing estimators.
Cramer-Rao-Frechet (CRF) introduced methods to calculate those
inequalities, and in few special cases they calculated the absolute
minimums of variances.
%
%
Reference~\cite{landi08} was
dedicated to extend the CRF approach to inequalities
for estimators of the least-squares methods for Gaussian
probability density functions (PDFs).
The Gaussian PDFs were used in a toy-model discussed in ref.~\cite{landi08,landi07}.
That model was defined in ref.~\cite{landi07}
to illustrate  the essential improvements introduced in the least-squares method
by the appropriate hit variances (heteroscedasticity), and to compare them with
the standard least squares without those differences (homoscedasticity).
However, it was clear to us, well before ref.~\cite{landi08,landi07},
that realistic effective variances of the observations (hits)
drastically improved the quality of track fitting, in  models of silicon trackers.
References~\cite{landi05,landi06} clearly showed the importance of
calculations dedicated to the estimations of the hit properties
relevant to track reconstruction.
Among these fitting results, one was very interesting~\cite{landi06}:
a linear growth of the maximums of the momentum distributions
with the number N of detecting layers. Such growth is much faster than the
corresponding growth, as $\sqrt{\mathrm{N}}$, of the standard least squares.
Further details were added in ref.~\cite{landi07} to the origin
of this linear growth and to show its generality in tracker physics.
For this task, it was simulated the estimator of the direction,
of a set of straight tracks, in a realistic model of silicon tracker.
To better illustrate these new results, a section
of ref.~\cite{landi07} was dedicated to a Gaussian toy-model
very easily reproducible, by the interested reader, with
few lines of MATLAB code~\cite{MATLAB}.
The Gaussian toy-model gives results very similar to the more complex
realistic simulations. Its linear growth has a very effective visual
impact for the unusual slenderness of the parameter distributions and
their fast growth with N. Instead, the standard least-squares does not move
from its modest growth as $\sqrt{\mathrm{N}}$.  Such evident differences
activated, in few of our readers (and referees), reminiscences of mathematical
statistics where the $\sqrt{\mathrm{N}}$ is a frequent result.
In particular, the widely reported form of the CRF-bound was used
to raise doubts~\cite{unk}  on the realizability of the linear growths,
in spite of a their reproduction of our toy-model.
Evidently, a simulation is a numerical demonstration, and it
rules out any theorem (if it could exist) with conflicting statements.
However, the necessity to remove forever any possible criticism about our
growth (linear or faster than linear) was the motivation of our ref.~\cite{landi08}.
Given that the criticisms were essentially addressed  to the Gaussian toy-model,
the demonstrations were based on Gaussian PDFs.
Reference~\cite{landi08} proved the perfect consistency of
those results with the CRF inequalities and the non optimality of the
standard least squares estimators for Gaussian PDFs.
The analytical properties of the Gaussian PDFs are ideal for those
developments, and, in general, the CRF method requires special
type of PDFs defined as regular models~\cite{mathstat,particle_group,econometrics}.
Even if the inequalities turn out to be independent from their
original PDFs, the extension to non-Gaussian PDFs, or irregular
models, can only be guessed by the analytical forms of the variances.
Given that a large part of our approaches are base on
truncated Cauchy-like PDFs (typical irregular models), the absence of  direct
demonstrations could allow the survival of possible doubts.
It is our aim to eliminate any misunderstanding
with a direct demonstration of the variance inequalities even for irregular models.
To illustrate the effects of these results, we introduce a
toy-model whose PDFs have  rectangular shapes.
This PDF is used in ref.~\cite{mathstat}, and elsewhere, as prototype of irregular
model. The parameters of the toy-model
are identical to those of ref.~\cite{landi08} with the sole differences in the
form of the hit PDFs.
The rectangular toy-model shows strong similarity with the Gaussian
one, as suggested by the formal identities of the new
inequalities with those for the Gaussian PDFs.

\begin{figure} [h!]
\begin{center}
\includegraphics[scale=0.75]{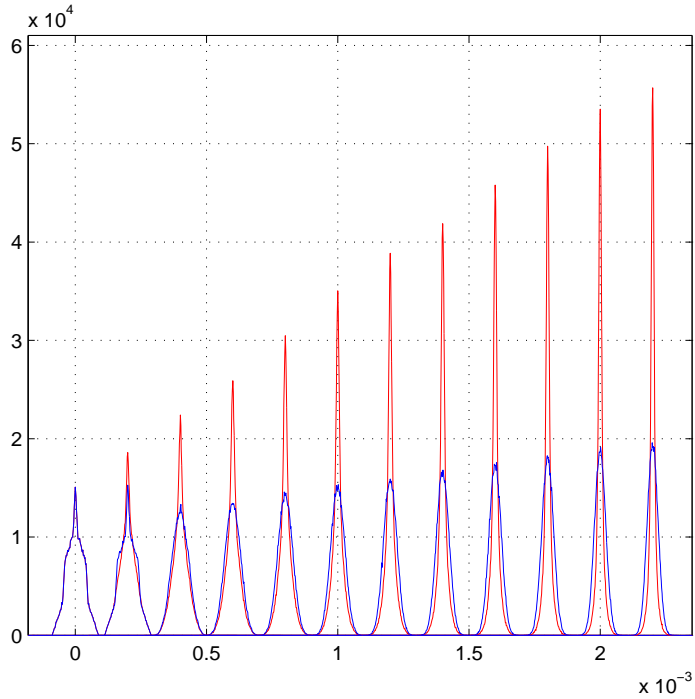}
\caption{\em Empirical PDFs of the fitted (straight)
track direction for tracker models with N=2 to N=13 detecting layers.
The PDFs for N=2 are centered on zero, the others are shifted
by N-2 steps of fixed amplitude. The blue PDFs are the
results of the standard least-squares. The red PDFs are
given by the weighted least-squares. (The reported dimensions are pure numbers)}
\label{fig:figure_1}
\end{center}
\end{figure}

Figure~\ref{fig:figure_1} shows  Monte Carlo simulations of a
simple tracker with N detecting layers of identical technology
(silicon microstrips) crossed at orthogonal
incidence by a set of parallel straight tracks of minimum
ionizing particles (MIPs). The reported estimator is the
track direction (the tangent of a small angle). The hit
uncertainties are obtained by a modulation of the
rectangular distributions ({\tt [rand-1/2]$\sqrt{12}$}), of
zero mean and unit variance, with two type of hit quality.
One of very good quality when the MIP crosses the detector
near to the borders ($20\%$ of strip width), and the other
of bad quality ($80\%$ of the strip width) around the strip
center. The standard deviation of the hits
are $\sigma_1=0.18$ (bad hits) and $\sigma_2=0.018$ strip
widths (the strip width is 63 $\mu\mathrm{m}$ and the
tracker length is 445 $\mathrm{mm}$). The hit properties are
a drastic simplification of the results of
refs.~\cite{landi05,landi06} for a type of silicon detector
of large use. (Results of ref.~\cite{CMSm} share
some similarity with these assumptions).  The simulated
orthogonal straight tracks
are fitted with two different least-squares methods:
one that uses only the hit positions (standard least
squares) and the second that uses also the hit standard
deviations (weighted least squares). Further details
of these simulations are reported in Appendix B of ref.~\cite{landi08}.
As it is evident in fig.~\ref{fig:figure_1}, the
blue lines  have a very slow growth as $\sqrt{\mathrm{N}}$.
Instead, the red lines of the heteroscedastic least-squares model
have an almost linear growth of the maximums of the empirical PDFs.
Translating these results on mean variances, the standard fit
shows a $1/\mathrm{N}$ trend, instead, the weighted fit shows a trend of the
mean variances as $1/\mathrm{N}^2$.
For precision, we have to specify better the origin of the linear growth, because
the inequalities only state relations among variances of
the estimators. This only implies that a set of
simulated tracks with heteroscedastic hit errors, reconstructed with the
appropriate estimators, have smaller mean variances
than the mean variance of the corresponding estimators of 
the standard least-squares. The linear growth with $N$
is mainly due to the usual simulation of a tracker
that distributes, on the tracks, good or bad hits with a binomial
probability distribution. In the weighted least squares, the track resolution
is mainly dominated by the good hits that, in the binomial distribution, grow
as $N$ in average (as far as the multiple scattering can be neglected).
Instead, the standard least-squares fit has a very weak sensitivity 
to the hit quality and to their binomial distributions. Its
variances have explicit $N$ dependencies (as $\sigma_m^2/{N}$
for the $\beta$ estimator) that modulate the maximums of the
distributions to grow as $\sqrt{N}$.

\section{Inequalities for irregular (and regular) models}

The developments of ref.~\cite{landi08} require very special
properties for the PDFs: continuity, at least two times
differentiability, inversion of derivation and integration, and
others. All these these properties are globally indicated
as "regular models". The Gaussian PDFs are perfect examples
of regular models. Instead, the rectangular PDFs (for example)
lack few of them and the CRF-methods can not be used in this
case. Even the schematic (realistic) models of
ref.~\cite{landi07,landi05,landi06} are irregular models
for the construction of the weights as described in
ref.~\cite{landi06}. Hence, a different method is required
to prove the variance inequalities for systems 
excluded by the CRF-assumptions.

\subsection{Definitions of the probability model}

Very general properties will be postulated for of the PDFs $f_j(x)$.
They cover regular and irregular models. To simplify the notations
of the integrations, the PDFs $f_j(x)$ are defined on all the
real line even for $f_j(x)$ different from zero on a finite interval.
Thus, the limits of the integrals are always $-\infty$ and $+\infty$,
these limits are ever subtended. The properties of the $f_j(x)$ are:
\begin{equation}\label{eq:equation_1}
\begin{aligned}
   & \ \theta\in R \ \ \ \ \ \ \ \ \  f_j(x-\theta)\geq 0 \ \ \ \ \ \ \ \ \ \ \ \
   \int f_j(x-\theta) \mathrm{d} x=1    \\
   &\int f_j(x-\theta)\, (x-\theta)\, \mathrm{d} x =0 \ \ \ \ \ \ \ \
   \int f_j(x-\theta)\, (x-\theta)^2\,\mathrm{d} x=\sigma_j^2
\end{aligned}
\end{equation}
The case of a single estimator is irrelevant for track
reconstruction and will be neglected. However, the variance
inequality can be directly proved by a simple Cauchy-Schwarz
inequality. Instead, we will concentrate first on straight tracks.
They cross N detection planes (microstrips)  producing the observations
$\{x_1,x_2\ldots x_N\}$ obtained by the PDFs $f_j(x_j-\beta-y_j\gamma)$.
Where $\beta$ is the track impact point, and $\gamma y_j$ is the
angular shift in a plane distant $y_j$ from the reference plane. 
The reference plane is selected to have $\sum_j\,y_j=0$.
The likelihood function for the N observations (hits) of the track is:
\begin{equation}\label{eq:equation_2}
    L_{1,2,\ldots,N}(\mathbf{x},\beta,\gamma)=f_1(x_1-\beta-y_1\gamma)
    f_2(x_2-\beta-y_2\gamma)\ldots\ldots
    f_N(x_N-\beta -y_N\gamma)\,.
\end{equation}
Two different mean-squares of observations are defined; the standard one
$M_s$ (homoscedastic) and the weighted one $M_w$ (heteroscedastic):
\begin{equation}\label{eq:equation_3}
    M_s=\sum_{j=1}^N (x_j-\beta-y_j\gamma)^2\ \  \  \  \
    \  \  \  \  \ M_w=\sum_{j=1}^N \frac{(x_j-\beta-y_j\gamma)^2}{\sigma_j^2} \,.
\end{equation}
The definition of the unbiased estimators is done with the derivatives
in $\beta$ and $\gamma$ of $M_s$ and $M_w$. The two required
vectors  are
$S_{1,2,\ldots,N}(\mathbf{x},\beta,\gamma)$ and
$U_{1,2,\ldots,N}(\mathbf{x},\beta,\gamma)$:
\begin{equation}\label{eq:equation_4}
    S_{1,2,\ldots,N}(\mathbf{x},\beta,\gamma)=-\frac{1}{2}
      \left (\begin{array} {l}
    \displaystyle{\frac{\partial M_s}{\partial\beta}}\\
    \displaystyle{\frac{\partial M_s}{\partial\gamma}}
     \end{array} \right )
     =
     \left (\begin{array} {l}
    \displaystyle{\sum_{j=1}^N {(x_j-\beta-y_j\gamma)}}\\
    \displaystyle{\sum_{j=1}^N {(x_j-\beta-y_j\gamma)y_j}}
     \end{array} \right ) \,.
\end{equation}
The condition $\sum_jy_j=0$ is neglected in $S_{1,2,\ldots,N}(\mathbf{x},\beta,\gamma)$
to maintain a similarity with $U_{1,2,\ldots,N}(\mathbf{x},\beta,\gamma)$.
\begin{equation}\label{eq:equation_5}
    U_{1,2,\ldots,N}(\mathbf{x},\beta,\gamma)=-\frac{1}{2}\left (\begin{array} {l}
    \displaystyle{\frac{\partial M_w}{\partial\beta}}\\
    \displaystyle{\frac{\partial M_w}{\partial\gamma}}
     \end{array} \right )
     =
     \left (\begin{array} {l}
    \displaystyle{\sum_{j=1}^N \frac{(x_j-\beta-y_j\gamma)}{\sigma_j^2}}\\
    \displaystyle{\sum_{j=1}^N \frac{(x_j-\beta-y_j\gamma)y_j}{\sigma_j^2}}
     \end{array} \right ) .
\end{equation}
The mean values of $S_{1,2,\ldots,N}(\mathbf{x},\beta,\gamma)$ and of
$U_{1,2,\ldots,N}(\mathbf{x},\beta,\gamma)$ with the
likelihood $L_{1,2,\ldots,N}(\mathbf{x},\beta,\gamma)$ are
always zero for eqs.~\ref{eq:equation_1}, and identically
for any of their linear combinations. The linear combinations
to extract the unbiased estimators from the vector
$S_{1,2,\ldots,N}(\mathbf{x},\beta,\gamma)$ and
$U_{1,2,\ldots,N}(\mathbf{x},\beta,\gamma)$ are given by
the following two $2\times 2$ matrices $\mathbf{R}$ and
$\mathbf{I}$ (from now on, the indications ${1,2,\ldots,N}$ of
the type of PDFs $f_j(x)$ will be subtended):
\begin{equation}\label{eq:equation_6}
\mathbf{R}=\frac{1}{2}\int \mathrm{d} \mathbf{x} L(\mathbf{x},\beta,\gamma)
       \left (\begin{array} {l l}
     \displaystyle{\frac{\partial^2 M_s}{\partial\beta^2}}
    &  \displaystyle{\frac{\partial^2 M_s}{\partial\beta\partial\gamma}}\\
     \displaystyle{\frac{\partial^2 M_s}{\partial\gamma\partial\beta}}
    & \displaystyle{\frac{\partial^2 M_s}{\partial\gamma^2}}
     \end{array}\right )
                        =\left (\begin{array} {l l}
                     \textstyle{N} & \ \ \ \textstyle{0}\\
                       &   \\
                      \textstyle{0} & \textstyle{\sum_{j=1}^N{y_j^2}}
                      \end{array} \right )\,, \ \ \ \
\end{equation}

\begin{equation}\label{eq:equation_7}
\mathbf{I}=\frac{1}{2}\int \mathrm{d} \mathbf{x} L(\mathbf{x},\beta,\gamma)
        \left (\begin{array} {l l}
     \displaystyle{\frac{\partial^2 M_w}{\partial\beta^2}}
    &  \displaystyle{\frac{\partial^2 M_w}{\partial\beta\partial\gamma}}\\
     \displaystyle{\frac{\partial^2 M_w}{\partial\gamma\partial\beta}}
    & \displaystyle{\frac{\partial^2 M_w}{\partial\gamma^2}}
        \end{array}\right )
                            =\left (\begin{array} {l l}
                     \textstyle{\sum_{j=1}^N\frac{1}{\sigma_j^2}} &  \textstyle{\sum_{j=1}^N\frac{y_j}{\sigma_j^2}}\\
                           &    \\
                      \textstyle{\sum_{j=1}^N\frac{y_j}{\sigma_j^2}} & \textstyle{\sum_{j=1}^N\frac{y_j^2}{\sigma_j^2}}
                          \end{array} \right )\,.
\end{equation}

\noindent
The matrices $\mathbf{R}$ and $\mathbf{I}$ are formally identical to the Fisher
information for Gaussian PDFs,  homoscedastic for $\mathbf{R}$ and heteroscedastic
for $\mathbf{I}$. The integral form of equation~\ref{eq:equation_7} serves to
save an analogy with ref.~\cite{landi08}, but it is irrelevant here.
The transposed matrices are indicated as $\mathbf{R}'$ or $\mathbf{I}'$.
For their symmetry, we often neglect to distinguish them from
their transposed.

The Gaussian model of ref.~\cite{landi08} shows that
equation~\ref{eq:equation_7} is even the integrals of the likelihood
$L(\mathbf{x},\beta,\gamma)$ with products $U_i\times {U_j}$.
This proof requires the the commutation of derivatives with integrations.
Unfortunately, these commutations are impossible for irregular models.
Instead, a direct calculation with the equations~\ref{eq:equation_1}
gives the last matrix of equation~\ref{eq:equation_7}.
The study of the variances requires the definition
of the unbiased estimators. The unbiased estimators
$T_1(\mathbf{x}),T_2(\mathbf{x})$ are obtained from
the vector $U$:
\begin{equation}\label{eq:equation_8}
  \textstyle{  U(\mathbf{x},\beta,\gamma)=\left (\begin{array} {c} \sum_{j=1}^N \frac{(x_j)}{\sigma_j^2} \\ \sum_{j=1}^N \frac{(x_j y_j)}{\sigma_j^2} \end{array} \right )-\mathbf{I}\left (\begin{array} {c} \beta \\ \gamma \end{array} \right ) \ \ \ \ \  \ \ \ \ \ \
    \mathbf{I}^{-1}U(\mathbf{x},\beta,\gamma)=\left (\begin{array} {c} T_1(\mathbf{x}) \\ T_2(\mathbf{x}) \end{array} \right )-\left (\begin{array} {c} \beta \\ \gamma \end{array} \right ) }
\end{equation}
The matrix $\mathbf{I}^{-1}$ and the two estimators $T_1(\mathbf{x}),T_2(\mathbf{x})$ are:
\begin{equation}\label{eq:equation_9}
    \mathbf{I^{-1}}=\frac{1}{\mathrm{Det}\,\{\mathbf{I}\}\,}\left (\begin{array} {l l}
                      \textstyle{\ \ \ \sum_{j=1}^N\frac{y_j^2}{\sigma_j^2}} &  \textstyle{-\sum_{j=1}^N\frac{y_j}{\sigma_j^2}}\\
                      \textstyle{-\sum_{j=1}^N\frac{y_j}{\sigma_j^2}} &  \textstyle{\ \ \ \sum_{j=1}^N\frac{1}{\sigma_j^2}} \end{array} \right ) \ \ \ \ \ \ \ \
    \left (\begin{array} {c} T_1(\mathbf{x}) \\ T_2(\mathbf{x}) \end{array} \right )=\mathbf{I}^{-1}\left (\begin{array} {c} \sum_{j=1}^N \frac{(x_j)}{\sigma_j^2} \\ \sum_{j=1}^N \frac{(x_j y_j)}{\sigma_j^2} \end{array} \right )
\end{equation}
For their special linear combinations of $U$ components, the integrals of $T_1(\mathbf{x})$ and $T_2(\mathbf{x})$ with the likelihood $L(\mathbf{x},\beta,\gamma)$ give $\beta$ and $\gamma$:
\begin{equation}\label{eq:equation_10}
\begin{aligned}
    &\int \big[T_{1}(\mathbf{x})-\beta\big]\,L(\mathbf{x},\beta,\gamma) \mathrm{d} \mathbf{x}=0 \ \ \ \ \ \ \ \ \ \ \ \ \ \ \int T_{1}(\mathbf{x})\,L(\mathbf{x},\beta,\gamma) \mathrm{d} \mathbf{x}=\beta\\
    &\int \big[T_{2}(\mathbf{x})-\gamma\big]\,L(\mathbf{x},\beta,\gamma) \mathrm{d} \mathbf{x}=0 \ \ \ \ \ \ \ \ \ \ \ \ \ \ \int T_{2}(\mathbf{x})\,L(\mathbf{x},\beta,\gamma) \mathrm{d} \mathbf{x}=\gamma\\
\end{aligned}
\end{equation}
As previously anticipated, the direct calculation of the integrals of the mean
values of the products $U_i\,U_j$ (with equations~\ref{eq:equation_1}) gives:
\begin{equation}\label{eq:equation_11}
\int \mathrm{d} \mathbf{x} L(\mathbf{x},\beta,\gamma)\cdot U(\mathbf{x},\beta,\gamma)\,U'(\mathbf{x},\beta,\gamma)
=\left (\begin{array} {l l}
                     \textstyle{\sum_{j=1}^N\frac{1}{\sigma_j^2}} &  \textstyle{\sum_{j=1}^N\frac{y_j}{\sigma_j^2}}\\
                      \textstyle{\sum_{j=1}^N\frac{y_j}{\sigma_j^2}} & \textstyle{\sum_{j=1}^N\frac{y_j^2}{\sigma_j^2}} \end{array} \right )\,,
\end{equation}
the matrix of equation~\ref{eq:equation_11} coincides with $\mathbf{I}$.
Differently from ref.~\cite{landi08}, no (impossible here) double differentiation
of the logarithm of the likelihood is required.
The variance-covariance matrix for $T_1(\mathbf{x}),T_2(\mathbf{x})$ becomes:
\begin{equation}\label{eq:equation_12}
    \int \mathbf{I^{-1}} U(\mathbf{x},\beta,\gamma)\cdot U'(\mathbf{x},\beta,\gamma)
    \mathbf{I^{-1}}'\,L(\mathbf{x},\beta,\gamma) \mathrm{d} \mathbf{x}
    =\mathbf{I^{-1}}'
\end{equation}
A very similar procedure is done on the vector $S(\mathbf{x},\beta,\gamma)$,
their unbiased estimators are:
\begin{equation}\label{eq:equation_13}
  \textstyle{  S(\mathbf{x},\beta,\gamma)=\left (\begin{array} {c} \sum_{j=1}^N {(x_j)} \\ \sum_{j=1}^N {(x_j y_j)} \end{array} \right )-\mathbf{R}\left (\begin{array} {c} \beta \\ \gamma \end{array} \right ) \ \ \ \ \  \ \ \ \ \ \
    \mathbf{R}^{-1}S(\mathbf{x},\beta,\gamma)=\left (\begin{array} {c} T_a(\mathbf{x}) \\ T_b(\mathbf{x}) \end{array} \right )-\left (\begin{array} {c} \beta \\ \gamma \end{array} \right ) }
\end{equation}
As previously done for other expressions, even the cross-correlation matrix of the estimators $T_1(\mathbf{x}),T_2(\mathbf{x})$
and $T_a(\mathbf{x}),T_b(\mathbf{x})$ is given by a direct integrations of the terms in square brackets:
\begin{equation}\label{eq:equation_14}
    \mathbf{R}^{-1}\big[\int S(\mathbf{x},\beta,\gamma)\cdot U'(\mathbf{x},\beta,\gamma)
    \,L(\mathbf{x},\beta,\gamma) \mathrm{d} \mathbf{x}\,\big]\,\mathbf{I^{-1}}'
    =\mathbf{R}^{-1} \left (\begin{array} {l l}
                     \textstyle{N} & \ \ \ \textstyle{0}\\
                      \textstyle{0} & \textstyle{\sum_{j=1}^N{y_j^2}} \end{array} \right ) \,\mathbf{I^{-1}}'
                      =\mathbf{I^{-1}}'\,.
\end{equation}
These results are identical to those obtained in ref.~\cite{landi08} for Gaussian
PDFs, but now they are derived with the sole use of the equations~\ref{eq:equation_1}, therefore valid for
irregular and regular models.

\subsection{The variance inequalities}
Here, we abandon the method of refs.~\cite{landi08,econometrics} and the use of
positive semi-definite matrices, for the more direct Cauchy-Schwarz inequality.
The matrix of equation~\ref{eq:equation_14} contains the integral:
\begin{equation}\label{eq:equation_15}
    \int \Big[T_a(\mathbf{x})-\beta\Big]\,\Big[T_1(\mathbf{X})-\beta\Big]
    \,L(\mathbf{x},\beta,\gamma) \mathrm{d} \mathbf{x}=(\mathbf{I}^{-1})_{1,1}\,,
\end{equation}
which implies the following  Cauchy-Schwarz inequality (with the
trivial substitution of the likelihood $L=\sqrt{L} \sqrt{L}$):
\begin{equation}\label{eq:equation_16}
    \int \Big[T_a(\mathbf{x})-\beta\Big]^2\,L(\mathbf{x},\beta,\gamma) \mathrm{d} \mathbf{x}
    \int \Big[T_1(\mathbf{X})-\beta\Big]^2
    \,L(\mathbf{x},\beta,\gamma) \mathrm{d} \mathbf{x}>(\mathbf{I}^{-1})_{1,1}^2\,.
\end{equation}
The equality is excluded because $T_a(\mathbf{x})$ and $T_1(\mathbf{x})$ are never
proportional for $N>2$ (and true heteroscedastic systems). For $N=2$, the minima
of $M_w$ and $M_s$ coincide  (we use this condition in fig.~\ref{fig:figure_1} to
check the consistency of the code).
Equation~\ref{eq:equation_12} states that the matrix element $(\mathbf{I}^{-1})_{1,1}$ is just
the variance of $T_1(\mathbf{x})$, therefore, the inequality becomes:
\begin{equation}\label{eq:equation_17}
\begin{aligned}
    &\int \Big[T_a(\mathbf{x})-\beta\Big]^2\,L(\mathbf{x},\beta,\gamma) \mathrm{d} \mathbf{x}\,\,
    (\mathbf{I}^{-1})_{1,1}>(\mathbf{I}^{-1})_{1,1}^2 \\
    &\int \Big[T_a(\mathbf{x})-\beta\Big]^2\,L(\mathbf{x},\beta,\gamma) \mathrm{d} \mathbf{x}>
    \int \Big[T_1(\mathbf{x})-\beta\Big]^2\,L(\mathbf{x},\beta,\gamma) \mathrm{d} \mathbf{x}\\
\end{aligned}
\end{equation}
The inequality for $T_b(\mathbf{X})$ and $T_2(\mathbf{x})$ is obtained with equation~\ref{eq:equation_15} to equation~\ref{eq:equation_17} and the matrix element $(\mathbf{I^{-1}})_{2\ 2}$, :
\begin{equation}\label{eq:equation_18}
\begin{aligned}
    &\int \Big[T_b(\mathbf{x})-\gamma\Big]^2\,L(\mathbf{x},\beta,\gamma) \mathrm{d} \mathbf{x}\,\,
    (\mathbf{I}^{-1})_{2,2}>(\mathbf{I}^{-1})_{2,2}^2 \\
    &\int \Big[T_b(\mathbf{x})-\gamma\Big]^2\,L(\mathbf{x},\beta,\gamma) \mathrm{d} \mathbf{x}>
    \int \Big[T_2(\mathbf{x})-\gamma\Big]^2\,L(\mathbf{x},\beta,\gamma) \mathrm{d} \mathbf{x}\,.\\
\end{aligned}
\end{equation}
Again, the equality is excluded because $T_a(\mathbf{x})$ and $T_1(\mathbf{x})$ are never
proportional for $N>2$.

Up to now, we supposed the absence of any effective variance in the standard least squares.
But often, the difference of detector technologies is introduced with a global effective
variance $\sigma_{a_l}$ for each hit in a detector layer.
It is easy to demonstrate that the precedent inequalities continue
to be valid even in this case. In fact, the correct likelihood
attributes to each hit its optimized variance, and this overcomes
the global effective variance. This
demonstration introduces small modifications with respect to
the previous one, essentially, the $\mathbf{R}$ matrix
contains now all these modifications. In any case, the
cross-correlation matrix is always generated by the product
$\mathbf{R}^{-1}\,\mathbf{R}\,\mathbf{I}^{-1}$.

\section{Extension of inequalities for momentum estimators to irregular models}

More important than the estimators for the straight
tracks are the estimators for fits of curved tracks in
magnetic field. Therefore, we will calculate these
inequalities for general irregular models. However, 
due to~\ref{eq:equation_1}, they remain valid even for 
regular models.

We will limit to a simpler case of high-momentum charged particles
in a tracker with homogeneous magnetic field orthogonal 
to the particle track. The segment 
of a circular path can be approximated with a parabola. 
As for straight tracks, the standard least-squares fit 
and the weighted one will be compared.
The likelihood function with the addition of curvature
$\eta$ for the N observations (hits) is:
\begin{equation}\label{eq:equation_19}
    L(\mathbf{x},\beta,\gamma\,,\eta)=f_1(x_1-\beta-y_1\gamma-y_j^2\eta)
    f_2(x_2-\beta-y_2\gamma-y_j^2\eta)\ldots\ldots
    f_N(x_N-\beta -y_N\gamma-y_N^2\eta)\,.
\end{equation}
The two mean squares of observations, $M_s$ (homoscedastic) and $M_w$
(heteroscedastic), are:
\begin{equation}\label{eq:equation_20}
    M_s=\sum_{j=1}^N (x_j-\beta-y_j\gamma -y_j^2\eta)^2\ \  \  \  \  \  \  \  \  \
    M_w=\sum_{j=1}^N \frac{(x_j-\beta-y_j\gamma-y_j^2\eta)^2}{\sigma_j^2}\,.
\end{equation}
The definition of the unbiased estimators requires the introduction
of the vectors $S(\mathbf{x},\beta,\gamma,\eta)$
and $U(\mathbf{x},\beta,\gamma,\eta)$:
\begin{equation}\label{eq:equation_21}
    S(\mathbf{x},\beta,\gamma\,.\eta)=-\frac{1}{2}
    \left (\begin{array} {l}
    \displaystyle{\frac{\partial M_s}{\partial\beta}}\\
    \displaystyle{\frac{\partial M_s}{\partial\gamma}}\\
    \displaystyle{\frac{\partial M_s}{\partial\eta}}
     \end{array} \right )
     =
     \left (\begin{array} {l}
    \displaystyle{\sum_{j=1}^N {(x_j-\beta-y_j\gamma-y_j^2\eta)}}\\
    \displaystyle{\sum_{j=1}^N {(x_j-\beta-y_j\gamma-y_j^2\eta)y_j}}\\
    \displaystyle{\sum_{j=1}^N {(x_j-\beta-y_j\gamma-y_j^2\eta)y_j^2}}
     \end{array} \right ) \,.
\end{equation}
To save the similarity with the form of $U(\mathbf{x},\beta,\gamma,\eta)$ the condition
$\sum_jy_j=0$ is not implemented. The vector $U(\mathbf{x},\beta,\gamma,\eta)$ is:
\begin{equation}\label{eq:equation_22}
    U(\mathbf{x},\beta,\gamma\,,\eta)=-\frac{1}{2}
     \left (\begin{array} {l}
    \displaystyle{\frac{\partial M_w}{\partial\beta}}\\
    \displaystyle{\frac{\partial M_w}{\partial\gamma}}\\
    \displaystyle{\frac{\partial M_w}{\partial\eta}}
     \end{array} \right )
     =
     \left (\begin{array} {l}
    \displaystyle{\sum_{j=1}^N \frac{(x_j-\beta-y_j\gamma-y_j^2\eta)}{\sigma_j^2}}\\
    \displaystyle{\sum_{j=1}^N \frac{(x_j-\beta-y_j\gamma-y_j^2\eta)y_j}{\sigma_j^2}}\\
    \displaystyle{\sum_{j=1}^N \frac{(x_j-\beta-y_j\gamma-y_j^2\eta)y_j^2}{\sigma_j^2}}
     \end{array} \right )\,.
\end{equation}
To extract the unbiased estimators the following two matrices are required:
\begin{equation}\label{eq:equation_23}
\mathbf{R}=\frac{1}{2}\int \mathrm{d} \mathbf{x} L(\mathbf{x},\beta,\gamma\,.\eta)
        \left (\begin{array} {l l l}
     \textstyle{\frac{\partial^2 M_s}{\partial\beta^2}}
    &  \textstyle{\frac{\partial^2 M_s}{\partial\beta\partial\gamma}}
    & \textstyle{\frac{\partial^2 M_s}{\partial\beta\partial\eta}}\\
     \textstyle{\frac{\partial^2 M_s}{\partial\gamma\partial\beta}}
    & \textstyle{\frac{\partial^2 M_s}{\partial\gamma^2}}
    & \textstyle{\frac{\partial^2 M_s}{\partial\gamma\partial\eta}}\\
     \textstyle{\frac{\partial^2 M_s}{\partial\eta\partial\beta}}
    & \textstyle{\frac{\partial^2 M_s}{\partial\eta\partial\gamma}}
    & \textstyle{\frac{\partial^2 M_s}{\partial\eta^2}}
       \end{array}\right )
                        =\left (\begin{array} {l l l}
                     \textstyle{N } &  \textstyle{0} & \textstyle{\sum_{j} y_j^2}\\
                      \textstyle{0} &  \sum_{j} y_j^2 & \sum_{j} y_j^3\\
                      \sum_{j} y_j^2 & \sum_{j} y_j^3 & \sum_{j} y_j^4
                          \end{array} \right )\,,
\end{equation}

\begin{equation}\label{eq:equation_24}
\mathbf{I}=\frac{1}{2}\int \mathrm{d} \mathbf{x} L(\mathbf{x},\beta,\gamma\,,\eta)
        \left (\begin{array} {l l l}
     \textstyle{\frac{\partial^2 M_w}{\partial\beta^2}}
    &  \textstyle{\frac{\partial^2 M_w}{\partial\beta\partial\gamma}}
    & \textstyle{\frac{\partial^2 M_w}{\partial\beta\partial\eta}}\\
     \textstyle{\frac{\partial^2 M_w}{\partial\gamma\partial\beta}}
    & \textstyle{\frac{\partial^2 M_w}{\partial\gamma^2}}
    & \textstyle{\frac{\partial^2 M_w}{\partial\gamma\partial\eta}}\\
     \textstyle{\frac{\partial^2 M_w}{\partial\eta\partial\beta}}
    & \textstyle{\frac{\partial^2 M_w}{\partial\eta\partial\gamma}}
    & \textstyle{\frac{\partial^2 M_w}{\partial\eta^2}}
        \end{array}\right )
                      =\left (\begin{array} {l l l}
                     \textstyle{\sum_j {1/\sigma_j^2}}  &  \textstyle{\sum_j{y_j/\sigma_j^2}} &
                     \textstyle{\sum_{j}{ y_j^2/\sigma_j^2}}\\
                     \textstyle{ \sum_j{y_j/\sigma_j^2}} & \textstyle{ \sum_{j}{y_j^2/\sigma_j^2}} & \textstyle{\sum_{j}{y_j^3/\sigma_j^2}}\\
                      \textstyle{\sum_{j}y_j^2/\sigma_j^2} & \textstyle{\sum_{j} y_j^3/\sigma_j^2} & \textstyle{\sum_{j} y_j^4/\sigma_j^2} \end{array} \right ).
\end{equation}

\noindent
Even here, the integrations with the likelihood are irrelevant. 
On these irregular models, the integrations
of products of $S\cdot S'$ and $U\cdot U'$ must be performed on 
their forms~\ref{eq:equation_21} and~\ref{eq:equation_22}.

The vectors $U$ and $S$ give
the unbiased estimators $T_1,T_2,T_3$ and $T_a,T_b,T_c$:
\begin{equation}\label{eq:equation_25}
  \textstyle{  U(\mathbf{x},\beta,\gamma,\eta)=\left (\begin{array} {c}
  \sum_{j=1}^N \frac{(x_j)}{\sigma_j^2} \\
  \sum_{j=1}^N \frac{(x_j y_j)}{\sigma_j^2} \\
   \sum_{j=1}^N \frac{(x_j y_j^2)}{\sigma_j^2}\end{array} \right )-
   \mathbf{I}\left (\begin{array} {c} \beta \\  \gamma \\ \eta \end{array} \right ) \ \ \ \ \  \ \
    \mathbf{I}^{-1}U(\mathbf{x},\beta,\gamma\eta)
    =\left (\begin{array} {c} T_1(\mathbf{x}) \\ T_2(\mathbf{x}) \\T_3(\mathbf{x}) \end{array} \right )-\left (\begin{array} {c} \beta \\ \gamma \\ \eta \end{array} \right ) } ,
\end{equation}
\begin{equation}\label{eq:equation_26}
  \textstyle{  S(\mathbf{x},\beta,\gamma,\eta)=
  \left (\begin{array} {c} \sum_{j=1}^N {(x_j)} \\
  \sum_{j=1}^N {(x_j y_j)}  \\
  \sum_{j=1}^N {(x_j y_j^2)}\end{array} \right )
  -\mathbf{R}\left (\begin{array} {c} \beta \\ \gamma \\ \eta \end{array} \right ) \ \ \ \ \  \ \
    \mathbf{R}^{-1}S(\mathbf{x},\beta,\gamma)=
    \left (\begin{array} {c} T_a(\mathbf{x}) \\ T_b(\mathbf{x}) \\  T_c(\mathbf{x})\end{array} \right )
    -\left (\begin{array} {c} \beta \\ \gamma \\ \eta \end{array} \right ) } .
\end{equation}
The direct integrations of the products $U_i \cdot U_j$ with the likelihood give:
\begin{equation}
    \int \mathrm{d} \mathbf{x} L(\mathbf{x},\beta,\gamma,\eta)\cdot U(\mathbf{x},\beta,\gamma,\eta)
    \,U'(\mathbf{x},\beta,\gamma,\eta)
= \mathbf{I}\,,
\end{equation}
and the variance-covariance matrix becomes:
\begin{equation}
    \mathbf{I^{-1}}\int \mathrm{d} \mathbf{x} L(\mathbf{x},\beta,\gamma,\eta) U(\mathbf{x},\beta,\gamma,\eta) \cdot
    \,U'(\mathbf{x},\beta,\gamma,\eta) \mathbf{I^{-1}}'
= \mathbf{I^{-1}}\,.
\end{equation}
The cross-correlation matrix of the estimators $T_1(\mathbf{x}),T_2(\mathbf{x}),T_3(\mathbf{x})$
and $T_a(\mathbf{x}),T_b(\mathbf{x}),T_c(\mathbf{x})$ is given by a direct calculations of the integrals.
The integrations in the square brackets give $\mathbf{R}$:
\begin{equation}\label{eq:equation_29}
    \mathbf{R}^{-1}\big[\int S(\mathbf{x},\beta,\gamma,\eta)\cdot U'(\mathbf{x},\beta,\gamma,\eta)
    \,L(\mathbf{x},\beta,\gamma,\eta) \mathrm{d} \mathbf{x}\,\big]\,\mathbf{I}^{-1}
=\mathbf{R}^{-1}\mathbf{R}\,\mathbf{I^{-1}}'=\mathbf{I}^{-1}
\end{equation}
and the cross-correlation of $T_c(\mathbf{x})$ with $T_3(\mathbf{x})$, the unbiased estimators of
the curvature, is:
\begin{equation}\label{eq:equation_30}
    \int \Big[T_c(\mathbf{x})-\eta\Big]\,\Big[T_3(\mathbf{X})-\eta\Big]
    \,L(\mathbf{x},\beta,\gamma,\eta) \mathrm{d} \mathbf{x}=(\mathbf{I}^{-1})_{3,3}\,.
\end{equation}
Equation~\ref{eq:equation_30} produces a Cauchy-Schwarz inequality:
(remembering that the likelihood can be written as $L=\sqrt{L} \sqrt{L}$):
\begin{equation}\label{eq:equation_31}
    \int \Big[T_c(\mathbf{x})-\eta\Big]^2\,L(\mathbf{x},\beta,\gamma,\eta) \mathrm{d} \mathbf{x}
    \int \Big[T_3(\mathbf{X})-\eta\Big]^2
    \,L(\mathbf{x},\beta,\gamma,\eta) \mathrm{d} \mathbf{x}>(\mathbf{I}^{-1})_{3,3}^2\,,
\end{equation}
the equality is excluded because $T_c(\mathbf{x})$ and $T_3(\mathbf{x})$ are never
proportional for $N>3$. For $N=3$, the minima of $M_w$ and $M_s$ coincide.
For equation~\ref{eq:equation_29}, the matrix element $(\mathbf{I}^{-1})_{3,3}$ is
the variance of $T_3(\mathbf{x})$. Therefore, the inequality for the curvature becomes:
\begin{equation}\label{eq:equation_32}
\begin{aligned}
    &\int \Big[T_c(\mathbf{x})-\eta\Big]^2\,L(\mathbf{x},\beta,\gamma,\eta) \mathrm{d} \mathbf{x}\,\,
    (\mathbf{I}^{-1})_{3,3}>(\mathbf{I}^{-1})_{3,3}^2 \\
    &\int \Big[T_c(\mathbf{x})-\eta\Big]^2\,L(\mathbf{x},\beta,\gamma\,,\eta) \mathrm{d} \mathbf{x}>
    \int \Big[T_3(\mathbf{x})-\eta\Big]^2\,L(\mathbf{x},\beta,\gamma\,,\eta) \mathrm{d} \mathbf{x}\\
\end{aligned}
\end{equation}
This inequality extends the validity of the results of our
simulations of ref.~\cite{landi06} to any heteroscedastic
system. The inequalities for $\beta$ and $\gamma$ can be
obtained with a similar procedure.

These demonstrations are directly extendible to a number
of parameters greater than those considered here.

\section{Effects of the inequalities on the resolution of the estimators}

Up to now we discussed the formal proofs of variance inequalities,
attributing to figure~\ref{fig:figure_1} the illustration of their
effects. We will show the direct connections of the variance inequalities
to the line-shapes of figure~\ref{fig:figure_1} and in the corresponding
figures of refs.~\cite{landi07,landi08}. To simplify the notations, we will
take $\beta$ and $\gamma$ equal to zero (as in the simulations).
It must be recalled that the illustrated simulations deal with a large
number of tracks with different successions of hit quality ($\sigma_{j_1},\cdots,\sigma_{j_N}$).
For each succession of $\sigma_{j_1},\cdots,\sigma_{j_N}$, extracted from the allowed set
of $\{\sigma_k\}$, we have a definite form of the estimators $T_1(\mathbf{x})$
and $T_2(\mathbf{x})$ (weighted least-squares). Instead, $T_a(\mathbf{x})$
and $T_b(\mathbf{x})$ (standard least-squares) are invariant. The studied estimator is
the track direction $\gamma$, thus $T_2(\mathbf{x})$ and $T_b(\mathbf{x})$ are
the selected estimators.

\subsection{The line-shapes of the estimators $T_2(\mathbf{x})$ and $T_b(\mathbf{x})$}

For any set of observations $x_{j_1},\cdots,x_{j_N}$,
the estimator $T_2(\mathbf{x})$gives a definite
value $\gamma^*$.
The set of all possible $\gamma^*$ has the probability distribution
$P_{j_1,\cdots,j_N}(\gamma^*)$ (with the method of ref.~\cite{landi06}):
\begin{equation}\label{eq:equation_33}
    P_{j_1,\cdots,j_N}(\gamma^*)=\int\, \mathrm{d} \mathbf{x}\, \delta\Big(\gamma^*-T_2(\sigma_{j_1},\cdots\sigma_{j_N},\mathbf{x})\Big)
    f_{j_1}(x_1)f_{j_2}(x_2)\cdots f_{j_N}(x_N)\,.
\end{equation}
This integral, for the linearity of $T_2$ in the observations $\mathbf{x}$,
gives a weighted convolution of the functions $f_{j_1}(x_1)f_{j_2}(x_2)\cdots f_{j_N}(x_N)$.
For the convolution theorems, the variance of $P_{j_1,\cdots,j_N}(\gamma^*)$
is given by the $(\mathbf{I^{-1}})_{2,2}$ of eq.~\ref{eq:equation_12} with
the corresponding $\sigma_j$:
\begin{equation}\label{eq:equation_34}
    \Sigma_{j_1,\cdots,j_N}=(\mathbf{I^{-1}})_{2,2}\,.
\end{equation}
In the case of ref.~\cite{landi08,landi07}, the convolution of Gaussian PDFs
gives another Gaussian function with the variance $\Sigma_{j_1,\cdots,j_N}$, 
therefore  $P_{j_1,\cdots,j_N}(\gamma^*)$ becomes:
\begin{equation}\label{eq:equation_35}
    P_{j_1,\cdots,j_N}(\gamma^*)=\exp\big[-\frac{{\gamma^*}^2}{2\, \Sigma_{j_1,\cdots,j_N}}\big]\frac{1}{\sqrt{2\,\pi\,\Sigma_{j_1,\cdots,j_N}}}\,.
\end{equation}
The mean value of $P_{j_1,\cdots,j_N}(\gamma^*)$ is zero, as it must be
for the unbiasedness of the estimator, and its maximum is ${1}/{\sqrt{2\,\pi\,\Sigma_{j_1,\cdots,j_N}}}$.

A similar procedure of eqs.~\ref{eq:equation_33}-~\ref{eq:equation_35} can be extended to the
standard least squares with the due differences $T_2\rightarrow T_b$,
$(\mathbf{I^{-1}})_{2,2}\rightarrow(\mathbf{R^{-1}})_{2,2}$ and $P_{j_1,\cdots,j_N}(\gamma^*)\rightarrow P_{j_1,\cdots,j_N}^b(\gamma_b^*)$.
In this case the variance of the convolution is defined:
\begin{equation}\label{eq:equation_36}
    S_{j_1,\cdots,j_N}=(\mathbf{R^{-1}})_{2,2}
\end{equation}
and, for Gaussian PDFs, the probability $P_{j_1,\cdots,j_N}^b(\gamma_b^*)$ is the function:
\begin{equation}\label{eq:equation_37}
    P_{j_1,\cdots,j_N}^b(\gamma_b^*)=\exp\big[-\frac{{\gamma_b^*}^2}{2\, S_{j_1,\cdots,j_N}}\big]\frac{1}{\sqrt{2\,\pi\,S_{j_1,\cdots,j_N}}}\,.
\end{equation}
The maximum of $P_{j_1,\cdots,j_N}^b(\gamma_b^*)$ is ${1}/{\sqrt{2\,\pi\,S_{j_1,\cdots,j_N}}} $.
For the inequality~\ref{eq:equation_18}, it is:
\begin{equation}\label{eq:equation_38}
    {1}/{\sqrt{2\,\pi\,S_{j_1,\cdots,j_N}}}\leq {1}/{\sqrt{2\,\pi\,\Sigma_{j_1,\cdots,j_N}}}\,.
\end{equation}
The equality is true if and only if all the $\sigma_j$ of the track are identical.

The law of large numbers allows to determine the convergence of the  line-shape toward the function:
\begin{equation}\label{eq:equation_39}
    \Pi(\gamma^*)=\sum_{j_1,\cdots,j_N}
     p_{j_1,\cdots,j_N}\,P_{j_1,\cdots,j_N}(\gamma^*) \ \ \ \ \ \ \ \
     \sum_{j_1,\cdots,j_N}\,p_{j_1,\cdots,j_N}=1,
\end{equation}
for the weighted fits. Analogously for the standard fits:
\begin{equation}\label{eq:equation_40}
    B(\gamma_b^*)=\sum_{j_1,\cdots,j_N}
     p_{j_1,\cdots,j_N}\,P_{j_1,\cdots,j_N}^b(\gamma_b^*)\,.
\end{equation}
The probability of the sequence of hit quality is $p_{j_1,\cdots,j_N}$.
The law of large numbers gives even an easy way for an analytical calculation
of the maximums of the two distributions:
\begin{equation}\label{eq:equation_41}
    \Pi(0)=\lim_{M\rightarrow \infty} \frac{1}{M} \sum_{k=1}^M {1}/{\sqrt{2\,\pi\,\Sigma_{k}}} \ \ \ \ \  \ \ \
    B(0)=\lim_{M\rightarrow \infty} \frac{1}{M} \sum_{k=1}^M {1}/{\sqrt{2\,\pi\,S_{k}}}\,,
\end{equation}
where the index $k$ indicates a track of a large number of simulated tracks,  $\Sigma_{k}$ is
the corresponding variance calculated as in the definition of $ (\mathbf{I^{-1}})_{2,2}$ and $S_{k}$
with the expression $ (\mathbf{R^{-1}})_{2,2}$. The inequality~\ref{eq:equation_38} assures
that $\ \Pi(0)\,>\,B(0)$.
For our 150000 tracks the convergence of the maximums of the empirical PDFs to
eqs.~\ref{eq:equation_41} is excellent for Gaussian PDFs. For the rectangular
PDFs of fig.~\ref{fig:figure_1}, the results of eqs.~\ref{eq:equation_41} are good for the
standard fit but slightly higher for weighted fits ($3\sim 4\%$). The strong similarity
of the line-shapes of the rectangular PDFs with the Gaussian PDFs is due to
the convolutions among various functions $f_j(x)$ that rapidly tend to approximate Gaussian
functions (Central Limit Theorem). The weighted convolutions of eq.~\ref{eq:equation_33}
has higher weights for the PDFs of good hits. These contribute mostly to the maximum of the
line-shape. Unfortunately, the good hits have a lower probability and the lower number
of the convolved PDFs produce a more approximate Gaussian than in the case of
standard fit. In fact, in the probability
$P_{j_1,\cdots,j_N}^b(\gamma_b^*)$ of the standard fit, the weights of the
convolved functions are identical and the convergence to a gaussian is better
(for N not too small).

\subsection{Resolution of estimators}

Figure~\ref{fig:figure_1} evidences the large differences of the direction PDFs
of the weighted least squares fits compared to the PDFs of the standard fits.
For $N>3$ the PDFs of the standard fits have negligible differences from a gaussian.
Instead, the PDFs of the weighted fits are surely non-Gaussian. The definition of
resolution becomes very complicated in these cases. In our previous papers, we
always used the maximums of the PDFs to compare different algorithms. Our convention was
lively contested by our readers of ref.~\cite{landi07}, they sustained that the usual
definition was  the standard deviation of all the data. It is evident the weakness of
this position for non Gaussian distributions. In fact, the standard deviation (as the
variance), for non Gaussian PDF, is principally controlled by
the tails of the distributions. Instead, the maximums are not directly effected
by the tails and they are true points of the PDFs.
In addition to this, a decreasing resolution, when the resolving power of the
algorithm/detector increases, easily creates contradictory statements.

\begin{figure} [h!]
\begin{center}
\includegraphics[scale=0.55]{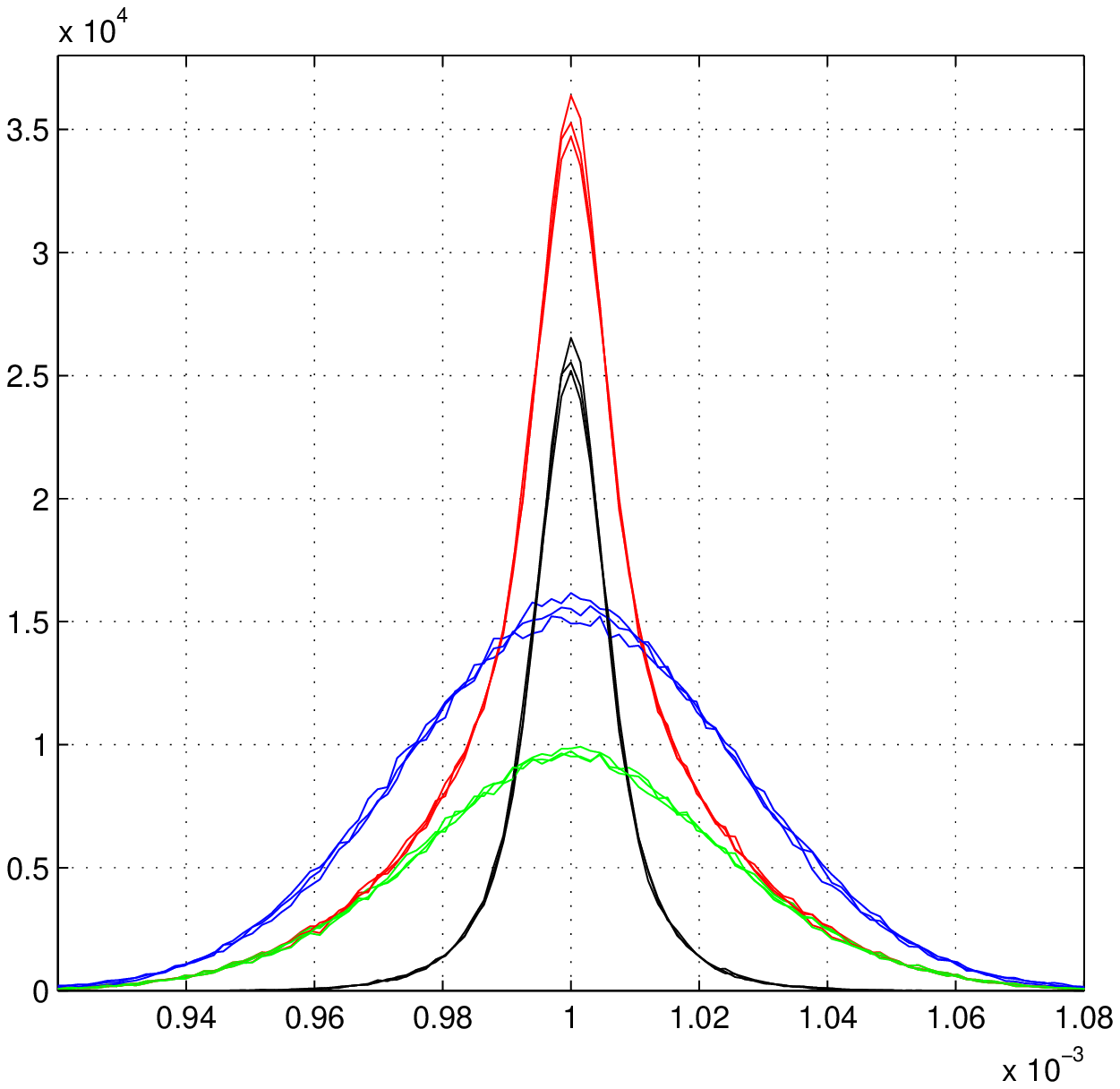}
\includegraphics[scale=0.58]{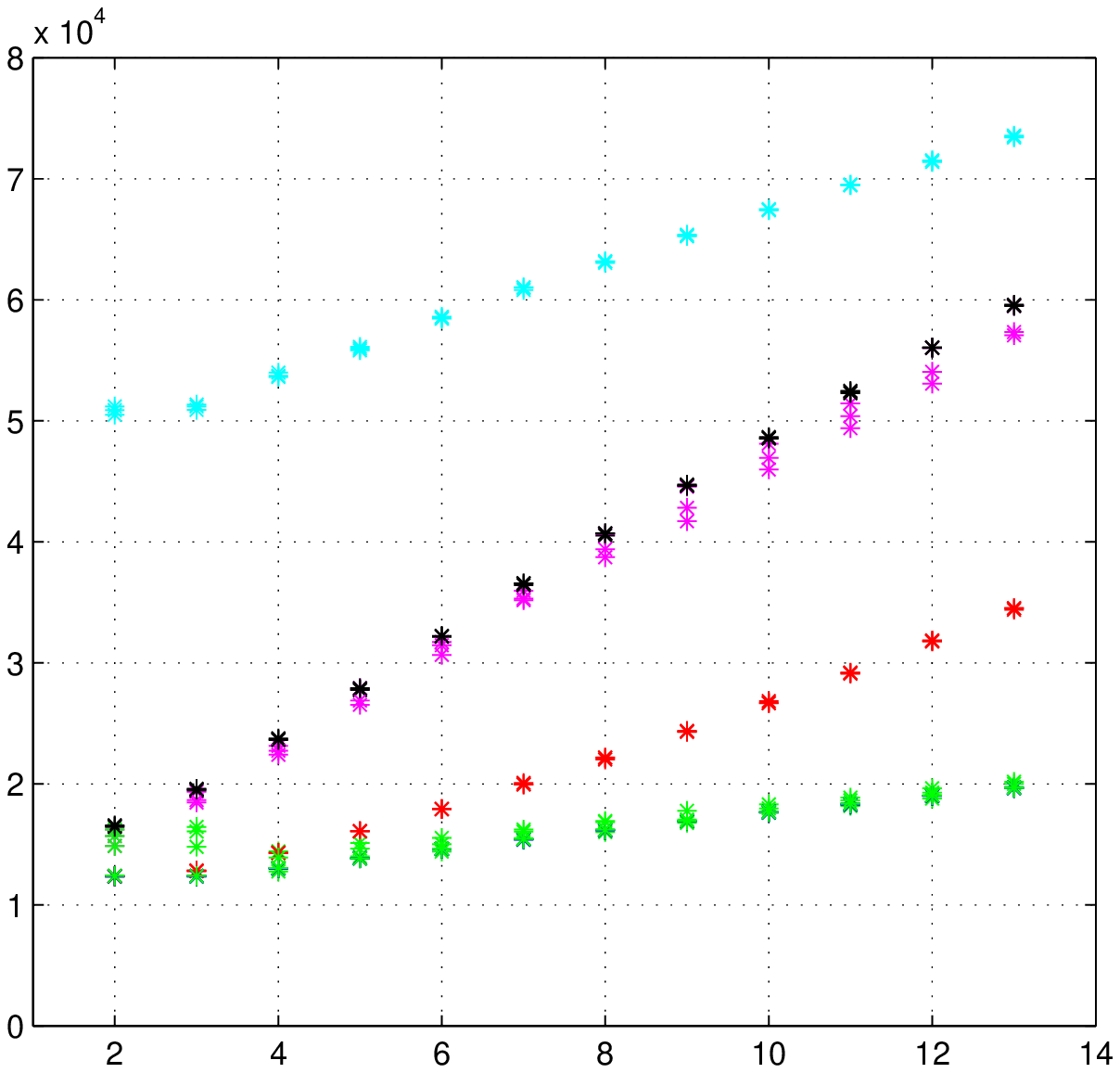}
\caption{\em Left plot: Empirical PDF of the fitted (straight)
track direction for a tracker model with N=7 detecting layers.
Red and blue lines as in fig.1. Black line, the fraction of red PDF due to tracks with
two or more good hits. Green line: the fraction of the red PDF due to one or zero good hits
The right plot reports the heights $1/\sqrt{2\,\pi\,\Sigma_l}$ where $\Sigma_l$ is a
variance for each N=$\{2,\cdots, 13\}$. The variances $\Sigma_l$ are:
green asterisks: variance of standard fits, red asterisks variances of weighted fits.
Ciano asterisks: the variance $\Sigma{eff}$.
Magenta asterisks maximums of red line of figure 1 (even for Gaussian
and triangular toy models), black asterisks: the results of $\Pi(0)$.
}
\label{fig:figure_2}
\end{center}
\end{figure}

Let us justify our preference.
The resolution indicates a property of an algorithm/instrument to discriminate
among near values of the calculated/measured parameters. It is evident that the
maximum discriminating power is obtained if the response of the algorithm/instrument
is a Dirac $\delta$-function. In this case, the parameters are confused if and only if
they coincide. Therefore, the algorithm with the best resolution is that
with the maximal similarity with a Dirac $\delta$-function.
Among a set of responses, it is evident that the
response with the higher maximum is the most similar to a Dirac $\delta$-function
and thus that with the largest resolution. This is the reason of our selection.
However, this is a conservative position. In fact, in recent years it is easy to find,
in literature,  more extremist conventions. It is a frequent practice (as in ref.~\cite{CMSb}
for example) to fit with a piece of Gaussian the central part (the core) of the distribution and
to take the standard deviation of this Gaussian as the measure of the resolution.\newline
The Gaussian toy-model allows to calculate the variance of a Gaussian
fitted to the core of the PDFs. A method to fit a Gaussian  is to fit
a parabola to the logarithm of the set of data. The second derivative of the parabola,
at the maximum of the distribution, is the inverse of the variance of the Gaussian.
This variance can be obtained by the function $\Pi(\gamma^*)$ of equation~\ref{eq:equation_39}
and equation~\ref{eq:equation_41}:
\begin{equation}\label{eq:equation_42}
   \frac{1}{\Sigma_{eff}}= \frac{d^2\ln\big[\Pi(\gamma^*)\big]}{d{\gamma^*}^2}\Big|_{\gamma^*=0}=
   \lim_{M\rightarrow \infty}
    \frac{ \sum_{k=1}^M {1}/{\sqrt{\Sigma_{k}^{3/2}}} }{ \sum_{n=1}^M {1}/{\sqrt{\Sigma_{n}}}}\,.
\end{equation}
The common factors are eliminated in the last fraction. Gaussian
functions, with standard deviations $\sqrt{\Sigma_{eff}}$, reproduce well
the core of $\Pi(\gamma^*)$ even if the Gaussian rapidly deviates from
$\Pi(\gamma^*)$. Instead, the Gaussian functions, with the height coinciding with
$\Pi(0)$, are outside the core of $\Pi(\gamma^*)$, their full width at half maximum
is always larger. Figure~\ref{fig:figure_2} illustrates few of these aspects for
three different toy-models.

The left side of fig~\ref{fig:figure_2} reports the line-shapes for N=7 layers
of the three toy models:
the Gaussian of refs.~\cite{landi08,landi07}, the rectangular of fig.~\ref{fig:figure_1}
and a toy-model with triangular PDFs. These three type of PDFs are those
considered by Gauss in his paper of 1823 on least squares method.
The line-shapes are
almost identical with small differences in their maximums. The red PDFs
are formed by the sums of the black and the green distributions, respectively
the fraction of the red PDF given by the tracks with two or more good hits,
and the fraction of the PDF given by tracks with
less than two good hits. It is evident the essential contribution
of the black distribution to the height of the red PDF. The PDFs with the lowest
maximums are those of the rectangular toy-model and the highest maximums
are those of the Gaussian toy-model.

The right side of  fig~\ref{fig:figure_2} reports the heights
$1/\sqrt{2\pi\Sigma_l}$ of  Gaussian PDFs with variance $\Sigma_l$
and the maximums of the empirical PDFs. The form of
equation~\ref{eq:equation_42} is a weighted mean of the inverse of the
variance of $P_{j_1,\cdots,j_N}(\gamma^*)$, with a weight proportional
to the maximum of its $P_{j_1,\cdots,j_N}(\gamma^*)$. This form allows
its extension to rectangular and triangular toy-models even if their
derivatives can be problematic. For the standard fits all these
heights are almost identical beyond N=3. The case N=2 is dealt in Sec.2.2, for N=3 the
standard fit tends to have the estimators very similar to those fof
N=2, due to $y_2=0$ that suppresses the central observations.
For the heteroscedastic fits all the effective heights $1/\sqrt{2\pi\Sigma_l}$
are largely different. Hence, the gain in resolution, with respect to the
homoscedastic fit, has values: 1.34, 2.4 and 4.1, for the model with seven
detecting layers. In this case, our preferred value is 2.4 times the 
result of the standard fit. However, we prefer to correlate the obtained 
PDFs to physical parameters of the problem. In ref.~\cite{landi06}, we used 
the magnetic field and the signal-to-noise ratio to measure the 
fitting improvements.\newline
At last, we have to remember that the non-Gaussian models allow a further 
increase in resolution with the maximum-likelihood search.
References~\cite{landi05,landi06} show the amplitude of the
increase in resolution due to the maximum-likelihood search.
If the triangular and rectangular-toy models have similar
improvements, the Gaussian-model would be the worst of all.

\section{Conclusions}

Inequalities among variances of the least-squares estimators are
calculated for heteroscedastic systems. They demonstrate that
the standard least squares method has not the minimal variances.
Thus, it is non optimal according to the usual definition of
optimality. These developments are explicitly constructed
for irregular models that can not
be handled with the Cramer-Rao-Frechet approach.
However, the conditions imposed to the probability
density functions are valid even for regular models. Thus, any
mixtures of regular and irregular models are equally dealt in these studies.
The linearity of the least-squares estimators (weighted or standard)
are key features for the demonstrations. To test the amplitudes
of the improvements due to these inequalities, the results of
simulations with one irregular models are reported. The simulation
has probability density functions of rectangular forms. This form
is used as an example of probability
model intractable with the Cramer-Rao-Frechet method.
The effects of the proved inequalities is connected to the formation
of the line-shape of the figures. The probability distributions for the
direction estimators is defined in the two different fitting methods.
The law of large numbers is used to reconstruct the maximums
of the plotted distributions. For the Gaussian distributions
the agreement if excellent. For the rectangular
distributions  slight differences are observed.
It is calculated even the Gaussian function
that fits the core of the model distributions.

%


\end{document}